\documentclass[11pt]{article}
\def\wtd{\widetilde}
\def\what{\widehat}
\def\qed{\vrule height7pt width5pt depth0pt}
\def\dfn{\stackrel {{\mbox{\scriptsize\rm{def}}}}{=}}

\def\bmx#1{\left(\begin{array}{@{}#1@{}}}
\def\emx{\end{array}\right)}
\newcommand{\brv}[1]{\setlength{\arraycolsep}{.5ex}
        \begin{array}{@{}#1@{}}(}
\newcommand{\erv}{)\end{array}}

\newtheorem{theorem}{Theorem}
\newtheorem{example}[theorem]{Example}

\usepackage{amssymb}
\usepackage{xcolor}
\usepackage{amsbsy}

\def\beps{\pmb{\epsilon}}

\def\br{\pmb{r}}
\def\bx{\pmb{x}}

\def\Blue{\textcolor{blue}}

\def\diag{{\rm diag}}
\def\subspan{{\rm span}}

\def\bbC{\mathbb{C}}

\def\cX{{\cal X}}

\def\sss{\scriptstyle}

\title{A Note on Eigenvalues of Perturbed Hermitian Matrices
    \footnote {\Blue{This paper originally published in {\em Linear Algebra and its Applications}, {\bf 395} (2005), 183--190.
               The main results are well received since.
               Unfortunately, the published version contains quite a number of typos. This arXiv version
               serves two purposes: 1) fix typos, and 2) add an appendix on how to best use the main results
               for sharper error bounds than straightforward applications in certain situations, e.g.,
               approximations from Krylov subspaces.}
               }}
\author{Chi-Kwong Li\thanks{Department of Mathematics,
          The College of William and Mary,
          Williamsburg, Virginia 23185
         ({\tt ckli@math.wm.edu}.)
         Research supported in part by the NSF grant DMS 0071994.}
\and Ren-Cang Li\thanks{Department of Mathematics,
        University of Kentucky,
        Lexington, KY 40506
        ({\tt rcli@ms.uky.edu}.)
        \Blue{
        Currently, Department of Mathematics,
        University of Texas at Arlington,
        Arlington, TX 76019
        ({\tt rcli@uta.edu}.)
        }
        This work was supported in part by
        the National Science Foundation CAREER award under Grant
        No. CCR-9875201.}
        }

\date{July 2004, \Blue{revised August 2025}}

\begin{document}

\maketitle

\begin{abstract}
Let
$$
A=\left(\begin{array}{cc} H_1 & E^*\\ E & H_2\end{array}\right) \quad
\hbox{ and } \quad
\wtd A=\left(\begin{array}{cc} H_1 & O\\ O & H_2\end{array}\right)$$
be two $N$-by-$N$ Hermitian matrices with eigenvalues $\lambda_1 \ge \cdots \ge
\lambda_{N}$ and
$\wtd \lambda_1 \ge \cdots \ge \wtd \lambda_N$, respectively.
Denote by $\|E\|$ the spectral norm of the matrix $E$, and
$\eta$ the spectral gap between the spectra of $H_1$ and $H_2$.
It is shown that
$$|\lambda_i - \wtd \lambda_i| \le
{2\|E\|^2 \over \eta+\sqrt{\eta^2+4\|E\|^2}} \, ,$$
which improves all the existing results.
Similar bounds are obtained for
singular values of matrices under block perturbations.
\end{abstract}



\section{Introduction}
\setcounter{equation}{0}
Consider a partitioned Hermitian matrix
\begin{equation}\label{eq:A}
A = \bordermatrix{       & \sss m & \sss n \cr
                  \sss m & H_1 & E^* \cr
                  \sss n & E & H_2 \cr },
\end{equation}
where $E^*$ is $E$'s complex conjugate transpose. At various situations
(typically when $E$ is {\em small}),
one is interested in knowing the impact of removing $E$ and $E^*$
on the eigenvalues of $A$. More specifically, one would like to obtain
bounds for the differences between that eigenvalues of $A$ and those of
its perturbed matrix
\begin{equation}\label{eq:tA}
\wtd A = \bordermatrix{  & \sss m & \sss n \cr
                  \sss m & H_1 & O  \cr
                  \sss n & O  & H_2 \cr }.
\end{equation}
Let $\lambda(X)$ be the spectrum of the square matrix $X$, and
let $\|Y\|$ be the spectral norm of a matrix $Y$, i.e., the
largest singular value of $Y$.
There are two kinds of bounds for the eigenvalues
$\lambda_1 \ge \cdots \ge \lambda_{m+n}$  and
$\wtd \lambda_1 \ge \cdots \ge \wtd \lambda_{m+n}$ of
$A$ and $\wtd A$, respectively:
\begin{enumerate}
    \item  \cite{bhat:96,parl:98,stsu:90}
     \begin{equation}\label{bd:classical-1}
     |\lambda_i- \wtd \lambda_i| \le \|E\|.
     \end{equation}
    \item       \cite{bhat:96,demm:97,govl:89,math:98,parl:98,stsu:90}
    If the spectra of $H_1$ and $H_2$ are disjoint, then
    \begin{equation}\label{bd:classical-2}
     |\lambda_i - \wtd \lambda_i| \le \|E\|^2/\eta,
     \end{equation}
     where
      $$
      \eta\dfn\min_{\mu_1\in\lambda(H_1),\,\mu_2\in\lambda(H_2)}|\mu_1-\mu_2|.
      $$
\end{enumerate}
The bound of the first kind does not use
information of the spectral distributions of the $H_1$ and $H_2$,
which will give (much) weaker bounds when
$\eta$ is not so small; while the bound of the second kind may blow up
whenever $H_1$ and $H_2$ have a common eigenvalue.
Thus both kinds have their own drawbacks, and it would be advantageous
to have bounds that are always
no bigger than $\|E\|$, of ${\cal O}(\|E\|)$ as $\eta\to 0$,
and at the same time behave like ${\cal O}(\|E\|^2/\eta)$ for not so
small $\eta$.
To further motivate our study, let us look at the following
$2\times 2$ example.

\smallskip
\begin{example}
Consider  the $2\times 2$ Hermitian matrix
\begin{equation}\label{eq:2by2H}
A=\left(\begin{array}{cc} \alpha & \epsilon \\
                          \epsilon & \beta\end{array}\right).
\end{equation}
\end{example}
Interesting cases are when $\epsilon$ is {\em small}, and thus
$\alpha$ and $\beta$ are {\em approximate\/} eigenvalues of
$A$. We shall analyze by how much the eigenvalues of $A$ differ from
$\alpha$ and $\beta$. Without loss of generality, assume
$$
\alpha>\beta,
\quad\mbox{and $\epsilon$ is real}.
$$
The eigenvalues of $A$, denoted by $\lambda_{\pm}$, satisfy
$\lambda^2-(\alpha+\beta)\lambda+\alpha\beta-\epsilon^2=0$; and thus
$$
\lambda_{\pm}=\frac {\alpha+\beta
                    \pm\sqrt{(\alpha+\beta)^2-4(\alpha\beta-\epsilon^2)}}2
      =\frac {\alpha+\beta
                    \pm\sqrt{(\alpha-\beta)^2+4\epsilon^2}}2.
$$
Now
\begin{eqnarray}
0<\left\{\begin{array}{c}\lambda_+-\alpha \\
                         \beta-\lambda_-\end{array}\right\}
  &=&\frac {-(\alpha-\beta)
                      +\sqrt{(\alpha-\beta)^2+4\epsilon^2}}2 \nonumber \\
  &=&\frac {2\epsilon^2}{(\alpha-\beta)
                       +\sqrt{(\alpha-\beta)^2+4\epsilon^2}}
            \label{eq:diff-2x2}
\end{eqnarray}
which provides a difference that enjoys the following properties:
$$
\frac {2\epsilon^2}{(\alpha-\beta)
                       +\sqrt{(\alpha-\beta)^2+4\epsilon^2}}
\left\{\begin{array}{l}
  \le \epsilon\quad\mbox{always,}\\
  \to \epsilon\quad\mbox{as $\alpha\to\beta^+$,}\\
  \le \epsilon^2/(\alpha-\beta).
  \end{array}\right.
$$
The purpose of this note is to extend this $2\times 2$ example and
obtain bounds which improve both (\ref{bd:classical-1})
and (\ref{bd:classical-2}).
Such results are not only
of theoretical interest but also important in the computations
of eigenvalues of Hermitian matrices \cite{kaha:66a,paig:74b,wilk:65}.

As an application, similar bounds are
presented for the singular value problem.

\section{Main Result} 
\setcounter{equation}{0}
\begin{theorem}\label{thm:mainThm}
Let
$$
A = \bordermatrix{  & \sss m & \sss n \cr
                  \sss m & H_1 & E^*  \cr
                  \sss n & E  & H_2 \cr } \quad \hbox{ and } \quad
\wtd A = \bordermatrix{  & \sss m & \sss n \cr
                  \sss m & H_1 & O  \cr
                  \sss n & O  & H_2 \cr }$$
be Hermitian matrices with eigenvalues
\begin{equation}\label{eq:eigvalsAtA}
\lambda_1\ge\lambda_2\ge\cdots\ge\lambda_{m+n}
\quad\mbox{and}\quad
\wtd\lambda_1\ge\wtd\lambda_2\ge\cdots\ge\wtd\lambda_{m+n},
\end{equation}
respectively.
Define
\begin{eqnarray}
  \eta_i &\dfn &\left\{\begin{array}{ll}
           \min\limits_{\mu_2\in\lambda(H_2)}|\wtd\lambda_i-\mu_2|,
               &\mbox{if $\wtd\lambda_i\in\lambda(H_1)$}, \vspace*{0.2cm}\\
           \min\limits_{\mu_1\in\lambda(H_1)}|\wtd\lambda_i-\mu_1|,
               &\mbox{if $\wtd\lambda_i\in\lambda(H_2)$},
        \end{array}\right. \label{eq:eta-i}\\
  \eta&\dfn&\min_{1\le i\le m+n}\eta_i
   =\min_{\mu_1\in\lambda(H_1),\,\mu_2\in\lambda(H_2)}|\mu_1-\mu_2|.
\label{eq:mineta-i}
\end{eqnarray}
Then for $i=1,2,\cdots, m+n$, we have
\begin{eqnarray}
|\lambda_i - \wtd \lambda_i|
  &\le&\frac {2\|E\|^2}{\eta_i+\sqrt{\eta_i^2+4\|E\|^2}}
        \label{ineq:mainIneq-1}\\
  &\le&\frac {2\|E\|^2}{\eta+\sqrt{\eta^2+4\|E\|^2}}.
        \label{ineq:mainIneq-2}
\end{eqnarray}
\end{theorem}
\it Proof. \rm
Suppose $U^*H_1U$ and $V^*H_2V$ are in the diagonal form with
their diagonal entries arranged in the descending order, respectively. We may assume
that $U = I_m$ and $V = I_n$. Otherwise, replace
$A$ by
$$
(U \oplus V)^* A (U \oplus V).
$$
We may perturb the diagonal of $A$ so that all entries are distinct,
and apply continuity argument for the general case.

\medskip
We prove the result by induction on $m+n$. If $m+n = 2$, the result
is clear (from our Example). Assume that $m+n > 2$, and the result is
true for Hermitian matrices of size $m+n-1$.

\medskip
First, refining an argument of Mathias \cite{math:98}, we show that
(\ref{ineq:mainIneq-1}) holds for $i = 1$.
Assume that the $(1,1)$st entry of $H_1$ equals $\wtd \lambda_1$.
By the min-max principle \cite{bhat:96,parl:98,stsu:90}, we have
$$
\lambda_1 \ge e_1^*Ae_1 = \wtd \lambda_1,
$$
where $e_1$ is the first column of the identity matrix.
No proof is necessary if $\lambda_1 = \wtd \lambda_1$. Assume
$\lambda_1 > \wtd \lambda_1$ and let
$$
X=\pmatrix{I_m & 0 \cr -(H_2-\lambda_1I_n)^{-1}E & I_n}.
$$
Then
$$
X^*(A-\lambda_1 I)X
= \left(\begin{array}{cc}
    M(\lambda_1) & 0 \\ 0 & H_2-\lambda_1 I_n\end{array}\right),$$
where
$$
M(\lambda) = H_1-\lambda I_m-E^*(H_2-\lambda I_n)^{-1} E.
$$
Since $A$ and $X^*AX$ have the same inertia,
we see that $M(\lambda_1)$
has zero as the largest eigenvalue.
Notice that the largest eigenvalue of $H_1 - \lambda_1 I$ is
$\wtd\lambda_1-\lambda_1\le 0$. Thus,
for $\delta_1 = |\lambda_1 - \wtd\lambda_1| = \lambda_1 - \wtd\lambda_1$,
we have (see  \cite[(10.9)]{parl:98})
$$
\lambda_1 \le \wtd \lambda_1 + \|E\|^2_2/(\delta_1 + \eta_1),
$$
and hence
$$\delta_1 \le
\|E\|^2/(\delta_1 + \eta_1).$$
Consequently,
$$
\delta_1 \le
{2\|E\|^2 \over \eta_1 + \sqrt{\eta_1^2 + 4\|E\|^2}}
$$
as asserted.
Similarly, we can prove the result if the $(1,1)$th entry of $H_2$ equals
$\wtd \lambda_1$. In this case, we will apply the inertia arguments to
$A$ and $YAY^*$ with
$$
Y = \pmatrix{I_m & 0 \cr -E(H_1 - \lambda_1 I_m)^{-1} & I_n}.
$$

\medskip
Applying the result of the last paragraph to $-A$, we see that
({\ref{ineq:mainIneq-1}) holds for $i = m+n$.

\medskip
Now, suppose $1 < i < m+n$. The result
trivially holds if $\lambda_i = \wtd \lambda_i$.
Suppose $\lambda_i \ne \wtd \lambda_i$.
We may assume that $\wtd \lambda_i > \lambda_i$. Otherwise,
replace $(A, \wtd A, i)$ by $(-A, -\wtd A, m+n-i+1)$.
Delete the row and column of $A$ that contain the diagonal entry
$\wtd \lambda_{m+n}$.
Suppose the resulted matrix $\what A$ has
eigenvalues $\nu_1 \ge \cdots \ge \nu_{m+n-1}$.
By the interlacing inequalities \cite[Section 10.1]{parl:98}, we have
\begin{equation}\label{eq:2}
\lambda_i \ge \nu_i \qquad \hbox{ and hence } \qquad
\wtd\lambda_i -\lambda_i \le \wtd \lambda_i - \nu_i.
\end{equation}
Note that $\wtd\lambda_i$ is the $i$th largest diagonal entries in
$\what A$.
Let $\what \eta_{i}$ be the minimum distance between $\wtd\lambda_i$
and the diagonal entries in the diagonal block $\what H_j$
in $\what A$ not containing $\wtd \lambda_i$,
where $j \in \{1,2\}$. Then
$$
{\what \eta}_{i} \ge \eta_i
$$
because  $\what H_j$ may have one fewer diagonal entries than $H_j$.
Let $\what E$ be the off-diagonal block of $\what A$. Then
$\|\what E\| \le \|E\|$.
Thus,
\begin{eqnarray*}
|\lambda_i - \wtd \lambda_i|
&=& \wtd \lambda_i - \lambda_i
    \hskip 0.999in  \hbox{ because } \wtd \lambda_i > \lambda_i \\
&\le& \wtd \lambda_i - \nu_{i}   \hskip 1.015in \hbox{ by } (\ref{eq:2}) \\
&\le&   {2\|\what E\|^2 \over \what \eta_{i} +
        \sqrt{\what\eta_{i}^2 + 4\|\what E\|^2}}
        \hskip 0.21in \hbox{ by induction assumption }\\
&\le&  {2\|\what E\|^2 \over \eta_{i} +
      \sqrt{\eta_{i}^2 + 4\|\what E\|^2}}
      \hskip 0.21in \hbox{ because } \what\eta_i \ge \eta_i \\
& = &  {1\over 2}\sqrt{\eta_i^2 + 4\|\what E\|^2} -\eta_i \\
& \le & {1\over 2} \sqrt{\eta_i^2 + 4\|E\|^2} -\eta_i
   \quad \hbox{ because } \|\what E \| \le \|E\| \\
&=& {2\|E\|^2 \over \eta_{i} + \sqrt{\eta_{i}^2 + 4\|E\|^2}}
\end{eqnarray*}
as asserted. \qed

\section{Application to Singular Value Problem}
\setcounter{equation}{0}

In this section, we apply the result in Section 2 to study singular
values of matrices.  For notational convenience in connection to our
discussion, we define the sequence of  singular
values of a complex $p \times q$ matrix $X$ by
$$\sigma(X) = (\sigma_1(X), \dots, \sigma_k(X)),$$
where $k = \max\{p,q\}$ and
$\sigma_1(X) \ge \cdots \ge \sigma_k(X)$ are the
nonnegative square roots of the eigenvalues of
the matrix $XX^*$ or $X^*X$  depending on which
one has a larger size. Note that the nonzero eigenvalues of
$XX^*$ and $X^*X$ are the same, and they give rise to the
nonzero singular values of $X$ which are of importance.
We have the following result concerning the nonzero
singular values of perturbed matrices.

\begin{theorem}\label{thm:mainThm-svd}
Let
$$
B = \bordermatrix{  & \sss k & \sss \ell \cr
                  \sss m & G_1 & E_1  \cr
                  \sss n & E_2  & G_2 \cr } \quad \hbox{ and } \quad
\wtd B = \bordermatrix{  & \sss k & \sss \ell \cr
                  \sss m & G_1 & O  \cr
                  \sss n & O  & G_2 \cr }
$$
be complex matrices with singular values
\begin{equation}\label{eq:singvalsBtB}
\sigma_1\ge\sigma_2\ge\cdots\ge\sigma_{\max\{m+n,k+\ell\}}
\quad\mbox{and}\quad
\wtd\sigma_1\ge\wtd\sigma_2\ge\cdots\ge\wtd\sigma_{\max\{m+n,k+\ell\}},
\end{equation}
respectively, so that $G_1$ and $G_2$ are non-trivial.
Define $\epsilon=\max\{\|E_1\|,\|E_2\|\}$, and
\begin{eqnarray}
  \eta_i &\dfn &\left\{\begin{array}{ll}
           \min\limits_{\mu_2\in\sigma(G_2)}|\wtd\sigma_i-\mu_2|,
               &\mbox{if $\wtd\sigma_i\in\sigma(G_1)$}, \vspace*{0.2cm}\\
           \min\limits_{\mu_1\in\sigma(G_1)}|\wtd\sigma_i-\mu_1|,
               &\mbox{if $\wtd\sigma_i\in\sigma(G_2)$},
        \end{array}\right. \label{eq:eta-i-sv}\\
  \eta&\dfn&\min_{1\le i\le \min\{m+n,k+\ell\}}\eta_i
       =\min_{\mu_1\in\sigma(G_1),\,\mu_2\in\sigma(G_2)}|\mu_1-\mu_2|.
      \label{eq:mineta-i-sv}
\end{eqnarray}
Then
for $i=1,2,\cdots, \min\{m+n,k+\ell\}$, we have
\begin{eqnarray}
|\sigma_i - \wtd \sigma_i|
  &\le&\frac {2\epsilon^2}{\eta_i+\sqrt{\eta_i^2+4\epsilon^2}}
        \label{ineq:mainIneq-1-sv}\\
  &\le&\frac {2\epsilon^2}{\eta+\sqrt{\eta^2+4\epsilon^2}} ,
        \label{ineq:mainIneq-2-sv}
\end{eqnarray}
and $\sigma_i = \wtd \sigma_i = 0$ for $i > \min\{m+n,k+\ell\}$.
\end{theorem}
{\sc Proof:}
By Jordan-Wielandt Theorem \cite[Theorem I.4.2]{stsu:90}, the eigenvalues
of
$$
\bmx{cc} O & B\\ B^* & O\emx
$$
are $\pm\sigma_i$ and possibly some zeros adding up to
$m+n+k+\ell$ eigenvalues. A similar statement holds for $\wtd B$.
Permuting the rows and columns appropriately, we see that
$$
\bmx{cc} O & B\\ B^* & O\emx
\quad\mbox{is similar to}\quad
\bmx{cc|cc} O     & G_1 & O     & E_1 \\
           G_1^* & O   & E_2^* & O   \\ \hline
           O     & E_2 & O     & G_2 \\
           E_1^*   & O   & G_2^* & O \emx,
$$
and
$$
\bmx{cc} O & \wtd B\\ \wtd B^* & O\emx
\quad\mbox{is similar to}\quad
\bmx{cc|cc} O     & G_1 &       &  \\
           G_1^* & O    &       &    \\ \hline
                 &      & O     & G_2 \\
                 &      & G_2^* & O \emx.
$$
Applying Theorem~\ref{thm:mainThm} with
$$
H_i=\bmx{cc} O     & G_i   \\
            G_i^* & O\emx  \quad \hbox{ and } \quad
E=\bmx{cc} O     & E_2 \\
           E_1^*   & O \emx ,
$$
we get the result. \qed

\bigskip
One can also apply the above proof to the degenerate cases when $G_1$
or $G_2$ in the matrix $B$ is trivial, i.e., one of the parameters
$m,n,k,\ell$ is zero. These cases are useful in applications.
We state one of them,  and one can easily extend it to other
cases.

\begin{theorem}
Suppose $B = \pmatrix{G \  E\cr}$ and
$\tilde B = \pmatrix{G \ O \cr}$ are
$p\times q$ matrices with
singular values
$$\sigma_1 \ge \cdots \ge \sigma_{\max\{p,q\}} \quad
\hbox{ and } \quad
\wtd \sigma_1 \ge \dots \ge \wtd\sigma_{\max\{p,q\}},$$
respectively.
Then for $i = 1, \dots, \min\{p,q\}$,
$$|\sigma_i - \wtd\sigma_i| \le
{2 \|E\|^2 \over 2\wtd\sigma_i + \sqrt{\wtd\sigma_i^2 + 4 \|E\|^2} }.$$
\end{theorem}

\section*{Acknowledgment} This paper evolved from an early note of
the second author in which Theorem~\ref{thm:mainThm} was proven for
special cases. The authors wish to thank Editor Volker Mehrmann for
making this collaboration possible.


\clearpage
\appendix
\section{\Blue{Discussions} (August 10, 2025)}
We have noticed that some recent articles on the similar topic use the main results of this paper in such a way that
the induced bounds are magnitude worse than the ``new'' bounds in the articles. It turns out that those discoveries
are misleading at best, if not anything else. In this appendix, we discuss how the main result in section~2 can
and should be used in those articles for  sharper error bounds than straightforward applications.
The same can be said about the result in section 3, but we omit the detail.

An equivalent form of Theorem 2 in terms of the eigen-residuals can be gotten as follows. Suppose that we have computed an $m$-dimensional approximate invariant subspace $\cX_1$ of Hermitian matrix $A\in\bbC^{N\times N}$ and let $X_1\in\bbC^{N\times m}$ be an orthonormal basis matrix of $\cX_1$, i.e., $\subspan(X_1)=\cX_1$, and $X_1^*X_1=I_m$. Let $H_1=X_1^*AX_1$, the Rayleigh quotient matrix
of $A$ with respect to  subspace $\cX_1$. As usual, we define the associated residual
as
\begin{equation}\label{eq:R}
R=AX_1-X_1H_1.
\end{equation}
Now expand $X_1$ to a unitary matrix $X:=[X_1,X_2]\in\bbC^{N\times N}$. We have
\begin{equation}\label{eq:XAX}
X^*AX=\left(\begin{array}{cc}
        X_1^*AX_1 & X_1^*AX_2 \\
        X_2^*AX_1 & X_2^*AX_2
      \end{array}\right)
      =:\left(\begin{array}{cc}
        H_1 & E^* \\
        E & H_2
      \end{array}\right)
\end{equation}
which is in the form of (\ref{eq:A}).  Theorem 2 now immediately applies to yield error bounds in terms of $\|E\|$.
We claim that
\begin{equation}\label{eq:E=R}
\|E\|=\|R\|
\end{equation}
and hence Theorem 2 also yields error bounds in terms of $\|R\|$.
To see (\ref{eq:E=R}), we notice that $\|R\|=\|X^*R\|=\|E\|$ because
$$
X^*R=X^*AX_1-X^*X_1H_1=\left(\begin{array}{c}
        X_1^*AX_1  \\
        X_2^*AX_1
      \end{array}\right)-\left(\begin{array}{c}
        H_1  \\
        0
      \end{array}\right)=\left(\begin{array}{c}
        0  \\
        E
      \end{array}\right).
$$

We observe that error bounds by Theorem 2 use the whole $E$ (or the whole $R$ for that matter) for all eigenvalues although  in (\ref{ineq:mainIneq-1}) they do use individual
gaps $\eta_i$.
It raises a question if somehow $E$ could be broken apart so that different pieces go with
different eigenvalues for sharper bounds. Indeed, that is possible and advantageous to do sometimes.
%
Without loss of generality, the basis matrix $X_1$ can be chosen so that
$$
H_1=\diag(\wtd\lambda_{j_1},\wtd\lambda_{j_2},\ldots, \wtd\lambda_{j_m}),
$$
in which case each column of
$$
X_1\equiv [\bx_1,\bx_2,\ldots,\bx_m]
$$
can be regarded as an approximate eigenvectors of $A$, and
$j_1<j_2<\cdots< j_m$ in terms of  earlier notation.
We still have (\ref{eq:XAX}). Define
\begin{equation}\label{eq:ri}
\br_i=A\bx_i-\bx_i\,\wtd\lambda_{j_i}\quad
\mbox{for $1\le i\le m$}.
\end{equation}
Introduce $\what H_i\in\bbC^{(N-1)\times (N-1)}$ standing for the matrix from striking out the $i$th row and column of
$X^*AX$ in (\ref{eq:XAX}) and write
$$
E=[\beps_1,\beps_2,\ldots,\beps_m].
$$
Theorem 2 says, for $1\le i\le m$,
\begin{equation}\label{eq:bdi}
|\lambda_{j_i}-\wtd\lambda_{j_i}|\le\frac {2\|\beps_i\|^2}{\what\eta_i+\sqrt{\what\eta_i^2+4\|\beps_i\|^2}}
    =\frac {2\|\br_i\|^2}{\what\eta_i+\sqrt{\what\eta_i^2+4\|\br_i\|^2}},
\end{equation}
 where
$$
\what\eta_i=\min_{\mu\in\lambda(\what H_i)} |\wtd\lambda_{j_i}-\mu|
\quad\mbox{for $1\le i\le m$}.
$$
Essentially, those individual error bounds are simply gotten by letting $m=1$ to begin with, for each eigenvalue of $H_1$.

The errors bounds by (\ref{eq:bdi}) use both individual gaps and different portions of $E$ (or different portions of $R$ for that matter). They
are particularly advantageous over the ones by Theorem 2 when $\cX_1$ is computed by
certain subspace methods for eigenvalue computations such as the Lanczos method, where
often those  $\|\beps_i\|=\|\br_i\|$ corresponding to a few largest and smallest eigenvalues of $H_1$ are {\em magnitude\/} smaller than the ones corresponding to other eigenvalues of $H_1$, yielding magnitude smaller error bounds
for the few largest and smallest eigenvalues, whereas error bounds straightforwardly from Theorem 2 may have similar orders
of magnitude
because of the use of the whole $E$ (or the whole $R$ for that matter).

How do the error bounds by Theorem 2 compare against those in (\ref{eq:bdi})? It is a tradeoff between mathematical elegance and practical sharpness, to put it simply. Yet, mathematically, (\ref{eq:bdi}) is no more general than
(\ref{ineq:mainIneq-1}). In fact, (\ref{ineq:mainIneq-1}) is the stronger one in mathematics because it implies
(\ref{eq:bdi}), even though the latter delivers much smaller error bounds in those articles referred to at the beginning of this appendix.

\end{document}